\documentclass[11pt]{article}

\usepackage{amsmath, amssymb, amsthm, latexsym}

\newcommand{\bbbn}{{\mathbb N}}
\newcommand{\words}{\{0,1\}^*}
\newcommand{\cantor}{\{0,1\}^\bbbn}
\newcommand{\segment}{\!\upharpoonright \!}

\theoremstyle{plain}
\newtheorem{theorem}{Theorem}
\newtheorem{proposition}[theorem]{Proposition}
\newtheorem{corollary}[theorem]{Corollary}

\theoremstyle{definition}
\newtheorem{definition}[theorem]{Definition}

\theoremstyle{remark}

\title{Is Randomness ``native" to Computer Science?}
\author{
Marie Ferbus-Zanda
\footnote{LIAFA, CNRS \& Universit\'e Paris 7,
2, pl. Jussieu 75251 Paris Cedex 05 France}\\
    \newcounter{fnnumber}
    \setcounter{fnnumber}{\value{footnote}}
{\footnotesize\tt ferbus@logique.jussieu.fr}
\and
Serge Grigorieff
{\footnotemark[\value{fnnumber}]}\\
{\footnotesize\tt seg@liafa.jussieu.fr}
}
\begin{document}
\date{August 2003}
\maketitle
\noindent
{\small
{\em 
%
Original paper published in}\\
Current Trends in Theoretical Computer Science.
\\
G. Paun, G. Rozenberg, A. Salomaa (eds.).
\\
World Scientific, Vol.2, p.141--180, 2004
\medskip\\
{\em Earlier version in}
\\
Yuri Gurevich's ``Logic In Computer Science Column"\\
Bulletin of EATCS, vol. 74, p.78--118, June 2001
}
{\footnotesize \textnormal \tableofcontents}
\newpage
%
\section{From probability theory to Kolmogorov complexity}
%
\subsection{Randomness and Probability theory}
\noindent {\bf Quisani}\footnote{Quisani is a student
with quite eclectic scientific curiosity,
who works under Yuri Gurevich's supervision }: \
I just found a surprising assertion
on Leonid Levin's home page:
\begin{quote}
{\em While fundamental in many areas of Science,
randomness is really ``native" to Computer Science.}
\end{quote}
\noindent
Common sense would rather consider randomness
as intrinsically relevant to Probability theory!

\medskip  \noindent {\bf Authors}: \
Levin also adds:
{\em ``The computational
nature of randomness was clarified by Kolmogorov."}

The point is that,
from its very origin to modern axiomatization
around 1933 \cite{kolmo33} by
Andrei Nikolaievitch Kolmogorov (1903--1987),
Probability theory carries a paradoxical result:
\begin{quote}\em
if we toss an unbiaised coin 100 times then
100 heads are just as probable as any other outcome!
\end{quote}
\noindent
As Peter G\'acs pleasingly remarks
(\cite{gacs93}, p. 3),
{\em
this convinces us only that the axioms of
Probability theory, as developped in \cite{kolmo33},
do not solve all mysteries
that they are sometimes supposed to.}
\medskip\\
In fact, since Laplace, much work has been devoted
to get {\em a mathematical theory of random objects},
notably by Richard von Mises (1883--1953)
(cf. \S\ref{vonmises}).
But none was satisfactory up to the 60's when such a theory
emerged on the basis of computability.
\\
As it sometimes occurs, the theory was discovered by several
authors independently.\footnote{For a detailed analysis
of {\em who did what, and when},
see \cite{li-vitanyi} p.89--92.}
In the USA, Ray J. Solomonoff (b. 1926),
1964 \cite{solomonoff} (a paper submitted in 1962)
and Gregory J. Chaitin (b. 1947), 1966
\cite{chaitin66}, 1969 \cite{chaitin69}
(both papers submitted in 1965).
In Russia, Kolmogorov, 1965 \cite{kolmo65},
with premisses announced in 1963 \cite{kolmo63}.

\medskip  \noindent {\bf Q}: \
Same phenomenon as for hyperbolic geometry with
Gauss, Lobatchevski and Bolyai.
I recently read a citation from Bolyai's father:
``When the time is ripe for certain things,
these things appear in different places in the manner of
violets coming to light in early spring".

\medskip  \noindent {\bf A}: \
Mathematics and poetry\ldots
Well, pioneered by
Kolmogorov, Martin-L\"of, Levin, G\'acs, Schnorr (in Europe)
and Chaitin, Solovay (in America),
the theory developped very fruitfully and
is now named {\em Kolmogorov complexity} or
{\em Algorithmic Information Theory}.

\medskip  \noindent {\bf Q}: \
So, Kolmogorov founded Probability Theory twice!
In the 30's and then in the 60's.

\medskip  \noindent {\bf A}: \
Hum\ldots
In the 30's Kolmogorov axiomatized Probability Theory
on the basis of measure theory,
i.e. integration theory on abstract spaces.
In the 60's, Kolmogorov
(and also Solomonoff and Chaitin independently)
founded a mathematical theory of {\em randomness}.
That it could be a new basis for Probability Theory
is not clear.

\medskip  \noindent {\bf Q}: \
What?
Randomness would not be {\em the} natural basis
for Probability Theory?

\medskip  \noindent {\bf A}: \
Random numbers are useful in different kinds of
applications:
simulations of natural phenomena,
sampling for testing ``typical case",
getting good source of data for algorithms,
\ldots (cf. Donald Knuth, \cite{knuth}, chapter 3).

\medskip \noindent
However, the notion of random object
as a mathematical notion is presently ignored
in lectures about Probability Theory.
Be it for the foundations or for the development of
Probability Theory,
such a notion is neither introduced nor used.
That's the way it is\ldots
There is a notion of random variable, but it
has really nothing to do with random objects.
Formally, they are just functions over some
probability space.
The name ``random variable" is a mere vocable
to convey the underlying {\em non formalized}
intuition of randomness.

\medskip  \noindent {\bf Q}: \
That's right.
I attended several courses on Probability Theory.
Never heard anything precise about random objects.
And, now that you tell me, I realize that there was
something strange for me with random variables.

\medskip
So, finally, our concrete experience of chance and
randomness on which we build so much intuition
is simply removed from the formalization of
Probability Theory.

\noindent
Hum\ldots Somehow, it's as if the theory of computability
and programming were omitting the notion of program,
real programs.

\noindent
By the way, isn' it the case?
In recursion theory, programs are reduced
to mere integers: G\"odel numbers!

\medskip  \noindent {\bf A}: \
Sure, recursion theory illuminates
but does not exhaust the subject of programming.

\medskip
As concerns a new foundation of Probability Theory,
it's already quite remarkable that Kolmogorov has
looked at his own work on the subject
with such a distance.
So much as to come to a new theory:
the mathematization of randomness.
However, it seems (to us) that Kolmogorov has been
ambiguous on the question of a new foundation.
Indeed, in his first paper on the subject
(1965, \cite{kolmo65}, p. 7),
Kolmogorov briefly evoked that possibility :
\begin{quote}\em
\dots to consider the use of the
[Algorithmic Information Theory] constructions
in providing a new basis for Probability Theory.
\end{quote}
However, later (1983, \cite{kolmo83}, p. 35--36),
he separated both topics
\begin{quote}\em
``there is no need whatsoever to change the
established construction of the mathematical
probability theory on the basis of the general theory
of measure.
I am not enclined to attribute the significance of
necessary foundations of probability theory to the
investigations [about Kolmogorov complexity] that
I am now going to survey.
But they are most interesting in themselves.
\end{quote}
though stressing the role of his new theory of random
objects for {\em mathematics as a whole}
(\cite{kolmo83}, p. 39):
\begin{quote}\em
The concepts of information theory as applied
to infinite sequences give rise to very interesting
investigations, which, without being indispensable
as a basis of probability theory, can acquire a
certain value in the investigation of the
algorithmic side of mathematics as a whole."
\end{quote}

  \noindent {\bf Q}: \
All this is really exciting.
Please, tell me about this approach to randomness.
%
\subsection{Intuition of finite random strings
and Berry's paradox}
\noindent {\bf A}: \
OK. We shall first consider finite strings.

\noindent
If you don't mind, we can start with
an approach which actually fails but conveys
the basic intuitive idea of randomness.
Well, just for a while, let's say that a finite
string $u$ is random if there is no shorter
way to describe $u$ but give the successive symbols
which constitute $u$.
Saying it otherwise, the shortest description
of $u$ is $u$ itself, i.e. the very writing
of the string $u$.

\medskip  \noindent {\bf Q}: \
Something to do with intensionality and
extensionality?

\medskip  \noindent {\bf A}: \
You are completely right.
Our tentative definition declares a finite string
to be random just in case it does not carry
{\em any} intensionality.
So that there is no description of $u$ but
the extensional one, which is $u$ itself.

\medskip  \noindent {\bf Q}: \
But the notion of description is
somewhat vague.
Is it possible to be more precise about
``description" and intensionality?

\medskip  \noindent {\bf A}: \
Diverse partial formalizations are possible.
For instance within any particular logical
first-order structure.
But they are quite far from exhausting the
intuitive notion of definability.
In fact, the untamed intuitive notion leads to
paradoxes, much as the intuition of truth.

\medskip  \noindent {\bf Q}: \
I presume you mean the liar paradox as concerns
truth.
As for definability, it should be Berry's paradox
about \\
\centerline{\em ``the smallest integer not definable
in less than eleven words"}
\noindent and this integer is indeed defined by this
very sentence containing only 10 words.

\medskip  \noindent {\bf A}: \
Yes, these ones precisely.
By the way, this last paradox was first mentioned
by Bertrand Russell, 1908 (\cite{russell}, p.222 or 150)
who credited G.G. Berry, an Oxford librarian,
for the suggestion.
%
\subsection{Kolmogorov complexity relative to a function}
\noindent {\bf Q}: \
And how can one get around such problems?

\medskip  \noindent {\bf A}: \
What Solomonoff, Kolmogorov and Chaitin did is a
very ingenious move: {\em instead
of looking for a general notion of definability,
they restricted it to computability}.
\noindent
Of course, computability is a priori as much a
vague and intuitive notion as is definability.
But, as you know, since the thirties, there is a
mathematization of the notion of computability.

\medskip  \noindent {\bf Q}: \
Thanks to Kurt, Alan and Alonzo.\footnote{Kurt G\"odel
                                           (1906--1978),
           Alan Mathison Turing (1912--1954),
           Alonzo Church (1903--1995).}

\medskip  \noindent {\bf A}: \
Hum\ldots
Well, with such a move,
general definitions of a string $u$ are
replaced by programs which compute $u$.

\medskip  \noindent {\bf Q}: \
Problem: we have to admit Church's thesis.

\medskip  \noindent {\bf A}: \
OK. In fact, even if Church's thesis were to
break down, the theory of computable functions
would still remain as elegant a theory as you
learned from Yuri and other people.
It would just be a formalization of a proper part
of computability, as is the theory of primitive
recursive functions or elementary functions.
As concerns Kolmogorov theory, it would still hold
and surely get extension to such a new context.

\medskip  \noindent {\bf Q}: \
But where do the programs come from?
Are you considering Turing machines or
some programming language?

\medskip  \noindent {\bf A}: \
Any partial computable function
$A:\{0,1\}^* \rightarrow \{0,1\}^*$
is considered as a programming language.
The domain of $A$ is seen as a family of programs,
the value $A(p)$ --- if there is any --- is the output
of program $p$.
As a whole, $A$ can be seen both
as a language to write programs
and as the associated operational semantics.

\noindent Now, Kolmogorov complexity relative
to $A$ is the function
$K_{A}:\{0,1\}^* \rightarrow \bbbn$
which maps a string $x$ to the length of shortest
programs which output $x$:
\begin{definition} \label{def:K}
$K_{A}(x) = \min\{ |p| \ :\ A(p)=x \}$
\end{definition}
\noindent
(Convention: $\min(\emptyset)=+\infty$,
so that $K_{A}(x)=+\infty$ if $x$ is outside
the range of $A$).

\medskip  \noindent {\bf Q}: \
This definition reminds me of a discussion I had
with Yuri some years ago (\cite{gurevich} p.76--78).
Yuri explained me things about Levin complexity.
I remember it involved time.

\medskip  \noindent {\bf A}: \
Yes. Levin complexity is a very clever variant of $K$
which adds to the length of the program
the log of the computation time to get the output.
It's a much finer notion. We shall not consider it
for our discussion about randomness.
You'll find some developments in \cite{li-vitanyi} \S 7.5.

\medskip  \noindent {\bf Q}: \
There are programs and outputs.
Where are the inputs?

\medskip  \noindent {\bf A}: \
We can do without inputs.
It's true that functions with no argument
are not considered in mathematics,
but in computer science, they are.
In fact, since Von Neumann, we all know that
there can be as much tradeof as desired
between input and program.
This is indeed the basic idea for universal
machines and computers.

\noindent
Nevertheless, Kolmogorov \cite{kolmo65} points a natural
role for inputs when considering
{\em conditional Kolmogorov complexity}
in a sense very much alike that of conditional
probabilities.

\noindent
To that purpose, consider a partial
computable function
$B:\{0,1\}^* \times \{0,1\}^* \rightarrow \{0,1\}^*$.
A pair $(p,y)$ in the domain of $B$ is interpreted as a
program $p$ together with an input $y$.
And $B(p,y)$ is the output of program $p$ on input $y$.
Kolmogorov \cite{kolmo65}
defines the conditional complexity relative
to $B$ as the function
$K_{B}(\ \mid\ ):\{0,1\}^* \times \{0,1\}^*
\rightarrow \bbbn$
which maps strings $x,y$ to the length of shortest
programs which output $x$ on input $y$:
\begin{definition} \label{def:condK}
$K_{B}(x\mid y)= \min \{ |p| \ :\ B(p,y)=x \}$
\end{definition}
%
\subsection{Why binary programs?}
\noindent {\bf Q}: \
Programs should be binary strings?

\medskip  \noindent {\bf A}: \
This is merely a reasonable restriction.

\noindent Binary strings surely have some flavor of machine level
programming.
But this has nothing to do with the present choice.
In fact, binary strings just allow for a fairness condition.
The reason is that Kolmogorov complexity deals with
lengthes of programs.
Squaring or cubing the alphabet divides all lengthes
by $2$ or $3$ as we see when going from binary
to octal or hexadecimal.
So that binary representation of programs is merely a
way to get an absolute measure of length.
If we were to consider programs $p$ written in some finite
alphabet $\Sigma$, we would have to replace
the length $|p|$ by the product
$|p|\log(card(\Sigma))$ where $card(\Sigma)$
is the number of symbols in $\Sigma$.
This is an important point when comparing
Kolmogorov complexities associated to diverse
programming languages, cf. \ref{sub:invariance}.
%
%
\subsection{What about other possible outputs?}
\noindent {\bf Q}: \
Outputs should also be binary strings?

\medskip  \noindent {\bf A}: \
Of course not.
In Kolmogorov approach, outputs are the finite objects
for which a notion of randomness is looked for.
Binary strings constitute a simple instance.
One can as well consider integers, or rationals or elements
of any structure $D$ with a natural notion of computability.
The modification is straightforward:
now $A$ is a partial computable function $A:\words\to D$ and
$K_A:D\to\bbbn$ is defined in the same way: $K_A(x)$, for $x\in D$,
is the minimum length of a program $p$ such that $f(p)=x$.
%
%
\section{Optimal Kolmogorov complexity}
%
\subsection{The Invariance Theorem} \label{sub:invariance}
\noindent {\bf Q}: \
Well, for each partial computable function $A:\words\to\words$
(or $A:\words\to D$, as you just explained)
there is an associated Kolmogorov complexity.
So, what is the Kolmogorov complexity of a given string?
It depends on the chosen $A$.

\medskip  \noindent {\bf A}: \
Now, comes the fundamental result of the theory, the so-called
{\em invariance theorem.}
We shall state it uniquely for Kolmogorov
complexity but it also holds for conditional
Kolmogorov complexity.

\noindent
Recall the {\em enumeration theorem}:
partial computable functions can be enumerated
in a partial computable way.
This means that there exists a partial computable function
$E:\bbbn\times\words\to\words$ such that,
for every partial computable function $A:\words\to\words$,
there is some $e\in\bbbn$ for which we have
$\forall p\ A(p)=E(e,p)$
(equality means that $A(p)$ and $E(e,p)$ are simultaneously
defined or not and, if defined, of course they must be equal).

\medskip  \noindent {\bf Q}: \
Wait, wait, I remember the diagonal argument which proves that
there is no enumeration of functions $\bbbn\to\bbbn$.
It goes through computability.
Given $E:\bbbn\times\words\to\words$, the function $A:\bbbn\to\bbbn$
such that $A(n)=E(n,n)+1$ is different from each one of the
function $n\mapsto E(n,e)$'s.
And if $E$ is computable then $A$ is computable.
\\
So, how can there be an enumeration of computable functions?

\medskip  \noindent {\bf A}: \
There is no computable enumeration of computable functions.
The diagonal argument you recalled proves that this is impossible.
No way.
\\
But, we are not considering computable functions
but partial computable functions.
This makes a big difference. The diagonal argument breaks down.
In fact, equality $E(n,n)=E(n,n)+1$ is not incoherent:
it just insures that $E(n,n)$ is not defined!

\medskip  \noindent {\bf Q}: \
Very strange property,indeed.

\medskip  \noindent {\bf A}: \
No, no. It's quite intuitive nowadays, in our world with computers.
Given a program in some fixed programming language, say language LISP,
an interpreter executes it.
Thus, with one more argument, the simulated program, the LISP compiler
enumerates all functions which can be computed by a LISP program.
Now, any partial computable function admits a LISP program.
Thus, the LISP interpreter gives you a computable enumeration of
computable functions.

\medskip  \noindent {\bf Q}: \
OK.

\medskip  \noindent {\bf A}: \
Let's go back to the invariance theorem.
We transform $E$ into a {\em one argument} partial computable
function $U:\words\to\words$ as follows:
Set
\medskip\\
$\left\{\begin{array}{rcll}
U(0^e1p)&=&E(e,p)
\\
U(q)&=&\mbox{undefined}&\mbox{if $q$ contains no occurrence of $1$}
\end{array}\right.$
\medskip\\
(where $0^e$ is the string $00...0$ with length $e$).
\\
Then, if $A:\words\to\words$ is partial computable
and $e_0$ is such that $\forall p\ A(p)=E(e_0,p)$, we have
\medskip\\
$\begin{array}{rclr}
K_{U}(x) & = & \min\{|q| \ :\ U(q)=x\}
& \mbox{(definition of }K_U)  \\
  & = &\min\{|0^e1p| \ :\ U(0^e1p)=x\}
&   \\
  & \leq &\min\{|0^{e_0}1p| \ :\ U(0^{e_0}1p)=x\}
& \mbox{(restriction to $e=e_0$)}  \\
  & = &\min\{|0^{e_0}1p| \ :\ A(p)=x\}
& \mbox{($e$ is a code for $A$)}  \\
  & = &e+1+\min\{|p| \ :\ A(p)=x\}
&   \\
  & = &e+1+ K_{A}(x)
& \mbox{(definition of }K_A \mbox{)}
\end{array}$
\medskip\\
Let's introduce useful notations.
For $f,g:\{0,1\}^* \rightarrow \bbbn$, let's write
$$
f\leq g +O(1)\ \ \ \ \mbox{(resp. $f = g +O(1)$)}
$$
to mean that there exists a constant c such that
$$
\forall x \ f(x)\leq g(x)+c\ \ \ \
\mbox{(resp. $\forall x \ |f(x)- g(x)|\leq c$)}
$$
i.e. $f$ is smaller than (resp. equal to)
$g$ up to an additive constant $c$.

\medskip
What we have just shown can be expressed as the following theorem
independently obtained  by
Kolmogorov (1965 \cite{kolmo65} p.5),
Chaitin (1966 \cite{chaitin66} \S9-11)
and Solomonoff (1964 \cite{solomonoff} p.12,
who gives the proof as an informal argument).
\begin{theorem}[Invariance theorem]
\label{thm:invariance}
$K_{U} \leq K_{A} +O(1)$
for any partial computable
$A:\{0,1\}^* \rightarrow\{0,1\}^*$.
In other words, up to an additive constant,
$K_{U}$ {\em is the smallest one} among the
$K_{A}$'s.
\end{theorem}
\noindent
Thus, up to an additive constant, there is a smallest $K_A$.
Of course, if $K_U$ and $K_V$ are both smallest,
up to an additive constant, then
$K_{U} = K_{V} + O(1)$.
Whence the following definition.
\begin{definition}[Kolmogorov complexity]
Kolmogorov complexity $K:\{0,1\}^* \rightarrow\bbbn$
is any fixed such smallest (p to an additive constant)
function $K_{U}$.
\end{definition}
Let's sum up.
The invariance theorem means that,
up to an additive constant,
there is an intrinsic notion of Kolmogorov complexity
and we can speak of {\em the}
Kolmogorov complexity of a binary string.

\medskip  \noindent {\bf Q}: \
Which is an integer defined up to a constant\ldots
Somewhat funny.

\medskip  \noindent {\bf A}: \
You witty!
Statements that only make sense in the limit
occur everywhere in mathematical contexts.

\medskip  \noindent {\bf Q}: \
Do not mind, I was joking.

\medskip  \noindent {\bf A}: \
In fact, Kolmogorov argued as follows
about the constant, \cite{kolmo65} p. 6:
\begin{quote}\em
Of course, one can avoid the indeterminacies
associated with the [above] constants,
by considering particular [\ldots functions $U$],
but it is doubtful that this can be done without
explicit arbitrariness.
One must, however, suppose that the different
``reasonable" [above universal functions] will lead
to ``complexity estimates" that will converge
on hundreds of bits instead of tens of thousands.
Hence, such quantities as the ``complexity" of
the text of ``War and Peace" can be assumed
to be defined with what amounts to uniqueness.
\end{quote}
\noindent {\bf Q}: \
Using the interpretation you mentioned a minute ago
with programming languages concerning the enumeration theorem,
the constant in the invariance theorem
can be viewed as the length of a LISP program which interprets $A$.

\medskip  \noindent {\bf A}: \
You are right.

%
\subsection{Coding pairs of strings}\label{coding}
\noindent {\bf A}: \
Have you noted the trick to encode an integer $e$ and a string $p$
into a string $0^e1p$ ?

\medskip  \noindent {\bf Q}: \
Yes, and the constant is the length of the extra part $0^e1$.
But you have encoded $e$ in unary. Why not use binary
representation and thus lower the constant?

\medskip  \noindent {\bf A}: \
There is a problem. It is not trivial to encode two binary strings
$u,v$ as a binary string $w$. We need a trick.
But first, let's be clear:
{\em encode} here means to apply a computable injective function
$\{0,1\}^* \times \{0,1\}^* \rightarrow\{0,1\}^*$.
\\
Observe that concatenation does not work: if $w=uv$
then we don't know which prefix of $w$ is $u$.
A new symbol $2$ inserted as a marker allows for
an encoding: from $w=u2v$ we can indeed recover $u,v$.
However, $w$ is no more a {\em binary} string.

\noindent
A simple solution uses this last idea together
with a padding function applied to $u$
which allows $1$ to become an end-marker.
Let $pad(u)$ be obtained by inserting
a new zero in front of every symbol in $u$.
For instance, $pad(01011)=0001000101$.
Now, a simple encoding $w$ of strings $u,v$
is the concatenation $w=pad(u)1v$.
In fact, the very definition of $pad$ insures that
the end of the prefix $pad(u)$ in $w$ is marked by
the first occurrence of $1$ at an odd position
(obvious convention: first symbol has position $1$).
Thus, from $w$ we get $pad(u)$ --- hence $u$ ---
and $v$ in a very simple way:
a finite automaton can do the job!
Observe that
\begin{equation}
|pad(u)1v| = 2|u|+|v|+1
\end{equation}
\medskip  \noindent {\bf Q}: \
Is the constant $2$ the best one can do?

\medskip  \noindent {\bf A}: \
No, one can iterate the trick.
Instead of padding $u$, one can pad the string
$\overline{|u|}$ which is the binary representation of
the length of $u$.
Look at the string $w=pad(\overline{|u|})1uv$.
The first occurrence of $1$ at an odd position
tells you which prefix of $w$ is $pad(\overline{|u|})$ and
which suffix is $uv$.
 From the prefix $pad(\overline{|u|})$,
you get $\overline{|u|}$, hence $|u|$.
 From $|u|$ and the suffix $uv$, you get $u$ and $v$.

\noindent
Nice trick, isn't it? And since
$|\overline{|u|}|= 1+\lfloor \log(|u|)\rfloor$,
we get
\begin{equation}
|pad(\overline{|u|})1uv| =
  |u|+|v|+2\lfloor \log(|u|)\rfloor +3
\end{equation}
\medskip  \noindent {\bf Q}: \
Exciting! One could altogether pad the length of
the length.

\medskip  \noindent {\bf A}: \
Sure.
$pad(\overline{|\overline{|u|}|})1\overline{|u|}uv$
is indeed an encoding of $u,v$.
The first occurrence of $1$ at an odd position
tells you which prefix of $w$ is
$pad(\overline{|\overline{|u|}|})$
and which suffix is $\overline{|u|}uv$.
 From the prefix $pad(\overline{|\overline{|u|}|})$ ,
you get $\overline{|\overline{|u|}|}$
hence $|\overline{|u|}|$.
Now, from $|\overline{|u|}|$
and the suffix $\overline{|u|}uv$
you get $\overline{|u|}$ --- hence $|u|$ --- and $uv$.
 From $|u|$ and $uv$, you get $u$ and $v$.
Also, a simple computation leads to
\begin{equation}
|pad(\overline{|\overline{|u|}|})1\overline{|u|}uv|
=
|u|+|v|+\lfloor \log(|u|)\rfloor +
2\lfloor\log(1+\lfloor \log(|u|)\rfloor)\rfloor +3
\end{equation}
\medskip  \noindent {\bf Q}: \
Of course, we can iterate this process.

\medskip  \noindent {\bf A}: \
Right.
But, let's leave such refinements.
%
\subsection{Non determinism}
\noindent {\bf Q}: \
Our problematic is about randomness.
Chance, randomness, arbitrariness,
unreasonned choice, non determinism
\ldots
Why not add randomness to programs
by making them non deterministic
with several possible outputs?

\medskip  \noindent {\bf A}: \
Caution: if a single program can output
every string then Kolmogorov complexity
collapses.
In order to get a non trivial theory, you
need to restrict non determinism.
There is a lot of reasonable ways to do so.
It happens that all lead to something
which is essentially usual Kolmogorov complexity
up to some change of scale
(\cite{shen-uspensky}, \cite{grigo-marion01}).
Same with the prefix Kolmogorov complexity
which we shall discuss later.
%
%
%
\section{How complex is Kolmogorov complexity?}
\noindent {\bf Q}: \
Well, let me tell you some points I see about $K_{A}$.

\noindent
The domain of $K_{A}$ appears to be the range
of $A$.
So that $K_{A}$ is total in case $A$ is onto.

\noindent
Since there are finitely many programs $p$ with length
$\leq n$, there can be only finitely many
$x$'s such that $K_{A}(x)\leq n$.
So that,
$\lim_{|x|\rightarrow +\infty} K_{A}(x)=+\infty$.

\medskip \noindent
Also, in the definition of $K_{A}(x)$,
there are 2 points:

\noindent
1) find some program which outputs $x$,

\noindent
2) make sure that all programs with shorter length
either do not halt or have output different from $x$.

\noindent
Point 2 does not match with definitions of
partial computable functions!
%
\subsection{Approximation from above}
\label{sub:approx}
\noindent {\bf A}: \
Right.
In general, $K_{A}$ {\em is not} partial computable.

\medskip  \noindent {\bf Q}: \
So, no way to compute $K_{A}$.

\medskip  \noindent {\bf A}: \
Definitely not, in general.
However, $K_{A}$ can be approximated from above:
\begin{proposition}
$K_{A}$ is the limit of a computable decreasing sequence
of functions.
\end{proposition}
\noindent
Moreover, we can take such a sequence of
functions with finite domains.

\noindent
To see this, fix an algorithm ${\cal A}$ for $A$
and denote $A_{t}$ the partial function obtained
by applying up to $t$ steps of algorithm ${\cal A}$
for the sole programs with length $\leq t$.
It is clear that $(t,p)\mapsto A_{t}(p)$ has
computable graph. Also, \\
\centerline{$K_{A_{t}}(x)=\min\{|p| \
:  \ p\in\{0,1\}^{\leq t} \mbox{ and } A_{t}(p)=x\}$}
so that $(t,x)\mapsto K_{A_{t}}(x)$ has computable graph too.

\noindent
To conclude, just observe that
$(t,x)\mapsto K_{A_{t}}(x)$ is decreasing in $t$
(with the obvious convention that $undefined = +\infty$)
and that
$K_{A}(x)=\lim_{t\rightarrow\infty} K_{A_{t}}(x)$.

\noindent
The same is true for conditional Kolmogorov
complexity $K_{B}(\ \mid\ )$.

\medskip  \noindent {\bf Q}: \
If $K_{A}$ is not computable, there should be no
computable modulus of convergence for
this approximation sequence.
So what can it be good for?

\medskip  \noindent {\bf A}: \
In general, if a function
$f: \{0,1\}^* \rightarrow \bbbn$
can be approximated from above by a computable sequence
of functions $(f_{t})_{t\in\bbbn}$ then
$X_{n}^f=\{x : f(x)\leq n \}$ is computably enumerable
(in fact, both properties are equivalent).
Which is a very useful property of $f$.
Indeed, such arguments are used in the proof of some
hard theorems in the subject of Kolmogorov complexity.

\medskip  \noindent {\bf Q}: \
Could you give me the flavor of what it can be
useful for?

\medskip  \noindent {\bf A}: \
Suppose you know that $X_{n}^f$ is finite
(which is indeed the case for $f=K_{A}$)
and has exactly $m$ elements then you can explicitly
get $X_{n}^f$.

\medskip  \noindent {\bf Q}: \
Explicitly get a finite set? what do you mean?

\medskip  \noindent {\bf A}: \
What we mean is that there is a computable
function which associates to any $m,n$ a code
(in whatever modelization of computability) for
a partial computable function which has range
$X_{n}^f$ in case $m$ is equal to the number
of elements in $X_{n}^f$.
This is not trivial. We do this thanks to
the $f_{t}$'s.

\noindent
Indeed, compute the $f_{t}(x)$'s for all
$t$'s and all $x$'s
until you get $m$ different
strings $x_{1},\ldots,x_{m}$ such that
$f_{t_{1}}(x_{1}),\ldots,f_{t_{m}}(x_{m})$
are defined and $\leq n$
for some $t_{1},\ldots,t_{m}$.

\noindent
That you will get such
$x_{1},\ldots,x_{m}$
is insured by the fact that $X_{n}^f$ has at least
$m$ elements and that
$f(x) =\min \{ f_{t}(x) \ : \ t\}$ for all $x$.

\noindent
Since $f\leq f_{t}$, surely these $x_{i}$'s
are in $X_{n}^f$.
Moreover, they indeed constitute the whole of $X_{n}^f$
since $X_{n}^f$ has exactly $m$ elements.
%
\subsection{Dovetailing}\label{sub:dovetailing}
\noindent {\bf Q}: \
You run infinitely many computations, some of which
never halt. How do you manage them?

\medskip  \noindent {\bf A}: \
This is called {\em dovetailing}.
You organize these computations
(which are infinitely many, some lasting forever)
as follows:\\
--- Do up to $1$ computation step of $f_{t}(x)$
for $0\leq t\leq 1$ and $0\leq |x|\leq 1$,\\
--- Do up to $2$ computation steps of $f_{t}(x)$
for $0\leq t\leq 2$ and $0\leq |x|\leq 2$,\\
--- Do up to $3$ computation steps of $f_{t}(x)$
for $0\leq t\leq 3$ and $0\leq |x|\leq 3$,\\
\ldots

\medskip  \noindent {\bf Q}: \
Somehow looks like Cantor's enumeration of $\bbbn^2$
as the sequence \\
\centerline{$(0,0)  \ \ \
(0,1) \ (1,0) \ \ \
(0,2) \ (1,1) \ (2,0)  \ \ \
(0,3) \ (1,2) \ (2,1) \ (3,0)\ldots $}

\medskip  \noindent {\bf A}: \
This is really the same idea. Here, it would rather be
an enumeration \`a la Cantor of $\bbbn^3$.
When dealing with a multi-indexed family of computations
$(\varphi_{t}(\vec x))_{\vec t}$ ,
you can imagine computation steps as tuples of integers
$(i,\vec t, \vec x)$
where $i$ denotes the rank of some
computation step of $f_{\vec t}(\vec x)$
(here, these tuples are triples).
Dovetailing is just a way of enumerating all points
in a discrete multidimensional space $\bbbn^k$
via some zigzagging \`a la Cantor.

\medskip  \noindent {\bf Q}: \
Well, Cantor wandering along a broken line which fills
the discrete plane
here becomes a way to sequentialize parallel
computations.
%
\subsection{Undecidability}\label{sub:undecidability}
\noindent {\bf Q}: \
Let's go back to $K_{A}$.
If $A$ is onto then $K_{A}$ is total.
So that if $A$ is also computable then $K_{A}$ should
be computable too.
In fact, to get $K_{A}(x)$ just compute
all $A(p)$'s, for increasing $|p|$'s
until some has value $x$:
this will happen since $A$ is onto.

\noindent
You said $K_{A}$ is in general undecidable.
Is this undecidability related to the fact that $A$
be partial computable and not computable?

\medskip  \noindent {\bf A}: \
No. It's possible that $K_{A}$ be quite trivial with
$A$ as complex as you want.
Let $f:\{ 0,1\}^* \rightarrow \bbbn$ be any
partial computable function.
Set $A(0x)=x$ and $A(1^{1+f(x)}0x)=x$.
Then, $A$ is as complex as $f$ though
$K_{A}$ is trivial since $K_{A}(x)= |x|+1$,
as is easy to check.

\medskip  \noindent {\bf Q}: \
Is it hard to prove that some $K_{A}$
is indeed not computable?

\medskip  \noindent {\bf A}: \
Not that much. If $U:\{0,1\}^* \rightarrow \{0,1\}^*$ is optimal
then we can show that $K_U$ is not computable.
Thus, $K$ (which is $K_U$ for some fixed optimal $U$) is not
computable.

\noindent
And this is where Berry's paradox comes back.
Consider the length-lexicographic order on binary strings:
$u<_{hier} v$ if and only if
$|u|<|v|$ or $|u|=|v|$ and $u$ is lexicographically before $v$.

Now, look, we come to the core of the argument.
The key idea is to introduce the function
$T:\bbbn \rightarrow \{0,1\}^*$ defined as follows:
$$
T(i)=
\mbox{ the } <_{hier} \mbox{ smallest } x \mbox{ such that }K_U(x)>i
$$
As you see, this function is nothing but an implementation
of the very statement in Berry's paradox
modified according to Kolmogorov's move from definability
to computability via the function $A$.
Clearly, we have
\begin{equation} \label{eq:1}
       K_U(T(i)) > i
\end{equation}
\noindent
Suppose, by way of contradiction, that $K_U$ is computable.
Then so is $T$ and so is the function $V:\words\to\words$ such that
$V(p)=T(Val_2(1p))$ where $Val_2(1p)$ is the integer with binary
representation $1p$.
\\
Now, if $i>0$ has binary representation $1p$ then $T(i)=V(p)$, so that
\begin{equation} \label{eq:2}
K_V(T(i))\leq |p|=\lfloor\log(i)\rfloor
\end{equation}
The invariance theorem insures that, for some $c$, we have
\begin{equation} \label{eq:3}
       K_U \leq K_V +c
\end{equation}
From inequalities (\ref{eq:1}), (\ref{eq:2}), (\ref{eq:3}) we get
$$
i < K_U(T(i)) \leq K_V(T(i))+c \leq \log(i) +c
$$
which is a contradiction for $i$ large enough since
$\lim_{i\rightarrow +\infty} \frac{\log(i)}{i} = 0$.
\\
Thus, our assumption that $K_U$ be computable is false.

\subsection{No non trivial computable lower bound}
\noindent {\bf Q}: \
Quite nice an argument.

\medskip  \noindent {\bf A}: \
Can get much more out of it:
\begin{theorem}
1) No restriction of $K_U$ to an infinite computable set
is computable.

\noindent
2) {\em Worse,}
if $X\subseteq \{0,1\}^*$ is computable and
$f: X\rightarrow \bbbn$ is a computable function
and $f(x) \leq  K_U(x)$ for all $x\in X$ then $f$ is bounded!
\end{theorem}
\noindent
To prove this, just change the above definition of
$T:\bbbn \rightarrow \{0,1\}^*$ as follows: \\
\centerline{$T(i)=
\mbox{ the } <_{hier} \mbox{ smallest } x\in X
\mbox{ such that } f(x)>i$}
Clearly, by definition, we have $f(T(i)) > i$.
Since $T(i)\in X$ and $f(x) \leq  K_U(x)$ for $x\in X$,
this implies equation (\ref{eq:1}) above.
Also, $f$ being computable, so are $T$ and $V$,
and equation (\ref{eq:2}) still holds.
As above, we conclude to a contradiction.

\medskip \noindent
Let's reformulate this result in terms of
the greatest monotonous
(with respect to $\leq_{hier}$) lower bound of $K_U$ which is

\centerline{
$m(x)= \min_{y\geq_{hier}x} K_U(y)$}
This function $m$ is monotonous and tends to $+\infty$
but it does so incredibly slowly:
{\em on any computable set it can not grow as fast
as any unbounded computable function.}

%
\subsection{Kolmogorov complexity and representation of objects}
\noindent {\bf Q}: \
You have considered integers and their base 2 representations.
Complexity of algorithms is often much dependent on the way objects
are represented. Here, you have not be very precise about
representation of integers.

\medskip  \noindent {\bf A}: \
There is a simple fact.
\begin{proposition}\label{prop:K-f(x)}
Let $f:\{0,1\}^* \rightarrow \{0,1\}^*$ be
partial computable.

\noindent 1) $K(f(x))\leq K(x) +O(1)$ for every $x$ in the domain of $f$.

\noindent 2) If $f$ is also injective then
$K(f(x)) = K(x) +O(1)$ for $x\in domain(f)$.
\end{proposition}
\noindent
Indeed, denote $U$ some fixed universal function
such that $K=K_{U}$.

\noindent
To get a program which outputs $f(x)$, we just
encode a program $\pi$ computing $f$
together with a program $p$ outputting $x$.

\noindent
Formally, let $A:\{0,1\}^* \rightarrow \{0,1\}^*$
be such that $A(pad(\pi)1z)$ is the output of
$\pi$ on input $U(z)$ for all $z\in \{0,1\}^*$.
Clearly, $A(pad(\pi)1p)=f(x)$
so that $K_{A}(f(x))\leq 2|\pi|+|p|+1$.
Taking $p$ such that $K(x)=|p|$, we get
$K_{A}(f(x))\leq K(x)+2|\pi|+1$.

\noindent
The Invariance Theorem insures that
$K(f(x))\leq K_{A}(f(x))+O(1)$, whence point 1 of the
Proposition.

\noindent
In case $f$ is injective, it has a partial computable
inverse $g$  with domain the range of $f$.
Applying point 1 to $f$ and $g$ we get
point 2.

\medskip  \noindent {\bf A}: \
So all representations of integers lead to the same Kolmogorov
complexity, up to a constant.

\medskip  \noindent {\bf A}: \
Yes, as long as one can computably go from one representation to the
other one.
%

%
\section{Algorithmic Information Theory}
%
\subsection{Zip/Unzip}
\noindent {\bf Q}: \
A moment ago, you said the subject was also
named Algorithmic Information Theory. Why?

\medskip  \noindent {\bf A}: \
Well, you can look at $K(x)$ as a measure of
the information contents that $x$ conveys.
The notion can also be vividly described
using our everyday use of
compression/decompression software
(cf. Alexander Shen's lecture \cite{shen}).
First, notice the following simple fact:
\begin{proposition} \label{prop:leqlength}
$K(x) \leq |x| +O(1)$
\end{proposition}
\noindent
Indeed, let $A(x)=x$. Then $K_{A}(x)=|x|$ and
the above inequality is a mere application
of the Invariance Theorem.

\medskip \noindent
Looking at the string $x$ as a file,
any program $p$ such that $U(p)=x$
can be seen as a {\em compressed file}
for $x$ (especially in case the right member
in Proposition \ref{prop:leqlength}) is indeed $<|x|$\ldots).

\noindent
So, $U$ appears as a {\em decompression algorithm}
which maps the compressed file $p$ onto the
original file $x$.
In this way, $K(x)$ measures the length of
the shortest compressed files for $x$.

\medskip \noindent
What does compression?
It eliminates redundancies,
explicits regularities to shorten the file.
Thus, maximum compression reduces the file to
the core of its information contents
which is therefore measured by $K(x)$.
%
\subsection{Some relations in
             Algorithmic Information Theory}
\noindent {\bf Q}: \
OK. And what does Algorithmic Information Theory look like?

\medskip  \noindent {\bf Q}: \
Conditional complexity should give some nice
relations as is the case with conditional
probability.

\medskip  \noindent {\bf A}: \
Yes, there are relations which have some
probability theory flavor.
However, there are often logarithmic extra terms
which come from the encoding of pairs of strings.
For instance, an easy relation:
\begin{eqnarray}\label{eq:cond}
K(x) \leq K(x\mid y) + K(y) +
2 \log(\min (K(x\mid y), K(y))) + O(1)
\end{eqnarray}
\noindent
The idea to get this relation is as follows.
Suppose you have a program $p$
(with no parameter) which outputs $y$
and a program $q$ (with one parameter)
which on input $y$ outputs $x$, then you can
mix them to get a (no parameter) program
which outputs $x$.

\noindent
Formally,
Suppose that $p,q$ are optimal, i.e.
$K(y)=|p|$ and $K(x \mid y)=|q|$.
Let $A_{1} , A_{2} :\{0,1\}^* \rightarrow \{0,1\}^*$
be such that\\
\centerline{
$A_{1}(pad(|z|)1zw)=A_{2}(pad(|w|)1zw)=V(w, U(z))$}
where $V$ denotes some fixed universal
function such that $K(\ \mid \ )=K_{V}(\ \mid \ )$.

\noindent
It is clear that
$A_{1}(pad(|p|)1pq)=A_{2}(pad(|q|)1pq)=x$,
so that\\
\centerline{
$K_{A_{1}}(x)\leq |p|+|q|+2\log(|p|)+O(1)$}
\centerline{
$K_{A_{2}}(x)\leq |p|+|q|+2\log(|q|)+O(1)$}
whence, $p,q$ being optimal programs,\\
\centerline{
$K_{A_{1}}(x)\leq K(y)+K(x \mid y)+2\log(K(y))+O(1)$}
\centerline{
$K_{A_{2}}(x)\leq K(y)+K(x \mid y)
+2\log(K(x \mid y))+O(1)$}
Applying the Invariance Theorem, we get (\ref{eq:cond}).
%
\subsection{Kolmogorov complexity of pairs}
\noindent {\bf Q}: \
What about pairs of strings in the vein of
the probability of a pair of events?

\medskip  \noindent {\bf A}: \
First, we have to define the Kolmogorov complexity
of pairs of strings.
The key fact is as follows:
\begin{proposition}\label{prop:pairs}
If $f,g :\{0,1\}^* \times \{0,1\}^* \rightarrow \{0,1\}^*$
are encodings of pairs of strings
(i.e. computable injections),
then $K(f(x,y))=K(g(x,y))+O(1)$.
\end{proposition}
As we always argue up to an additive constant,
this leads to:
\begin{definition}\label{eq:pairs}
The Kolmogorov complexity of pairs is
$K(x,y)=K(f(x,y))$ where $f$ is any fixed encoding.
\end{definition}
\noindent
To prove Proposition \ref{prop:pairs},
observe that $f\circ g^{-1}$ is a partial
computable injection such that $f=(f\circ g^{-1})\circ g$.
Then, apply Proposition \ref{prop:K-f(x)}
with argument $g(x,y)$ and function $f\circ g^{-1}$.

\subsection{Symmetry of information}
\noindent {\bf A}: \
Relation (\ref{eq:cond}) can be easily improved to
\begin{eqnarray}\label{eq:cond2}
K(x,y) \leq K(x\mid y) + K(y) +
2 \log(\min (K(x\mid y), K(y))) + O(1)
\end{eqnarray}
The same proof works. Just observe that
from both programs $pad(|p|)1pq$ and $pad(|q|)1pq$
one gets $q$ hence also $y$.

\medskip \noindent
Now, (\ref{eq:cond2}) can be considerably improved:
\begin{theorem}\label{thm:hard}
$|K(x,y)- K(x\mid y) - K(y)| = O(\log(K(x,y))$
\end{theorem}
\noindent
This is a hard result,
independently obtained by Kolmogorov and Levin
around 1967
(\cite{kolmo68}, \cite{zvonkin-levin} p.117).
We better skip the proof (you can get it in
\cite{shen} p. 6--7 or
\cite{li-vitanyi} Thm 2.8.2 p. 182--183).

\medskip  \noindent {\bf Q}: \
I don't really see the meaning of that theorem.

\medskip  \noindent {\bf A}: \
Let's restate it in another form.
\begin{definition}\label{def:mutual}
$I(x:y) = K(y) - K(y\mid x)$ is called the
algorithmic information about $y$ contained in $x$.
\end{definition}
\noindent
This notion is quite intuitive: you take the difference
between the whole information contents of $y$
and that when $x$ is known for free.

\noindent
Contrarily to what was expected in analogy
with Shannon's classical information theory,
this is not a symmetric function
However, up to a logarithmic term, it is symmetric:
\begin{corollary}\label{coro:mutual}
$|I(x:y) - I(y:x)| = O(log(K(x,y))$
\end{corollary}
\noindent
For a proof, just apply Theorem \ref{thm:hard} with
$K(x,y)$ and $K(y,x)$ and observe that
$K(x,y)=K(y,x)+O(1)$ (use Proposition \ref{prop:pairs}).
%
%
\section{Kolmogorov complexity and Logic}
%
\subsection{What to do with paradoxes}
\noindent {\bf Q}: \
Somehow, Solomonoff, Kolmogorov and Chaitin have built up
a theory from a paradox.

\medskip  \noindent {\bf A}: \
Right. In fact, there seems to be two
mathematical ways towards paradoxes.
The most natural one is to get rid of
them by building secured and delimited
mathematical frameworks which will leave
them all out (at least, we hope so\ldots).
Historically, this was the way followed
in all sciences.
A second way, which came up in the 20th century,
somehow integrates paradoxes into scientific
theories via some clever and sound (!) use of
the ideas they convey.
Kolmogorov complexity is such a remarkable
integration of Berry's paradox into mathematics.

\medskip  \noindent {\bf Q}: \
As G\"odel did with the liar paradox which
is underlying his incompleteness theorems.
Can we compare these paradoxes?

\medskip  \noindent {\bf A}:
Hard question.
The liar paradox is about truth while Berry's
is about definability.
Viewed in computational terms,
truth and definability somehow correspond to
denotational and operational semantics.

\noindent
This leads to expect connections between
incompleteness theorems \`a la G\"odel
and Kolmogorov investigations.
%
\subsection{Chaitin Incompleteness results }
\noindent {\bf Q}: \
So, incompleteness theorems can be obtained from
Kolmogorov theory?

\medskip  \noindent {\bf A}: \
Yes. Gregory Chaitin, 1971 \cite{chaitin71},
pointed and popularized a simple
but clever and spectacular application
of Kolmogorov complexity
(this original paper by Chaitin did not consider
$K$ but the number of states of Turing machines,
which is much similar).

\noindent
Let ${\cal T}$ be a computable theory containing
Peano arithmetic such that all axioms of
${\cal T}$ are true statements.
\begin{theorem}
There exists a constant $c$ such that if
${\cal T}$ proves $K(x)\geq n$ then $n\leq c$.
\end{theorem}
The proof is by way of contradiction and is a redo
of the undecidability of $K_{U_{PASCAL^{bin}}}$.
Suppose that ${\cal T}$ can prove statements
$K(x)\geq n$ for arbitrarily large $n$'s.
Consider a computable enumeration of all theorems
of ${\cal T}$ and let $f:\bbbn \rightarrow \bbbn$
be such that $f(n)$ is the first string $x$ such
that $K(x)\geq n$ appears as a theorem of ${\cal T}$.
Our hypothesis insures that $f$ is total, hence a
computable function.
By very definition,
\begin{equation}\label{eq:000}
K(\overline{f(n)}\geq n
\end{equation}
Also, applying Propositions \ref{prop:K-f(x)} and
\ref{prop:leqlength} we get
\begin{equation}\label{eq:001}
K(\overline{f(n)}) \leq K(\overline{n})+O(1)
\leq \log(n)+O(1)
\end{equation}
whence $n \leq \log(n)+O(1)$,
which is a contradiction if $n$ is large enough.

\medskip  \noindent {\bf Q}: \
Quite nice.
But this does not give any explicit statement.
How to compute the constant $c$ ?
How to get any explicit $x$'s such that $K(x)>c$ ?

\medskip  \noindent {\bf A}: \
Right. Hum\ldots you could also see this as a
particularly strong form of incompleteness: you have
a very simple infinite family of statements, only
finitely many can be proved
but you don't know which ones.
%
\subsection{Logical complexity of $K$}
\noindent {\bf A}: \
By the way, there is a point we should
mention as concerns the logical complexity
of Kolmogorov complexity.

\noindent
Since $K$ is total and not computable, its graph can not be
computably enumerable (c.e.).
However, the graph of any $K_{A}$ (hence that of $K$)
is always of the form $R \cap S$ where
$R$ is c.e. and $S$ is co-c.e.
(i.e. the complement of an c.e. relation).

\noindent
We can see this as follows.
Fix an algorithm {\cal P} for $A$ and denote $A^t$
the partial function obtained by applying up to $t$
computation steps of this algorithm. Then\\
\centerline{$K_{A}(x)\leq n \Leftrightarrow \exists t \
(\exists p\in \{0,1\}^{\leq n} \ A^t(p)=x)$}
The relation within parentheses is computable in $t,n,x$,
so that $K_{A}(x)\leq n$ is c.e. in $n,x$.

\noindent
Replacing $n$ by $n-1$ and going to negations,
we see that $K_{A}(x)\geq n$ is co-c.e.
Since
$K_{A}(x)=n \Leftrightarrow
(K_{A}(x)\leq n) \wedge (K_{A}(x)\geq n)$,
we conclude that $K_{A}(x)=n$ is the intersection of
an c.e. and a co-c.e relations.

\noindent
In terms of Post's  hierarchy, the graph of $K_{A}$ is
$\Sigma^0_{1}\wedge \Pi^0_{1}$, hence $\Delta^0_{2}$.
The same with $K_{B}(\ \mid\ )$.

\medskip  \noindent {\bf Q}: \
Would you remind me about Post's hierarchy?

\medskip  \noindent {\bf A}: \
Emil Post introduced families of relations
$R(x_{1}\ldots x_{m})$ on strings and/or integers.
Let's look at the first two levels:

$\Sigma^0_{1}$ and $\Pi^0_{1}$ are the respective families
of c.e. and co-c.e. relations,

$\Sigma^0_{2}$ is the family of projections of $\Pi^0_{1}$
relations,

$\Pi^0_{2}$ consist of complements of
$\Sigma^0_{2}$ relations.

\noindent
Notations $\Sigma^0_{i}$ and $\Pi^0_{i}$ come from the
following logical characterizations:

$R(\vec x)$ is $\Sigma^0_{1}$
if
$\ R(\vec x) \Leftrightarrow
\exists t_{1}\ldots\exists t_{k} \ T(\vec t, \vec x)\ $
with $T$ computable.

$R(\vec x)$ is $\Sigma^0_{2}$ if
$\ R(\vec x) \Leftrightarrow
\exists \vec t \ \forall \vec u\ T(\vec t, \vec u, \vec x)\ $
with $T$ computable.

$\Pi^0_{1}$ and $\Pi^0_{2}$ are defined similarly with
quantifications $\forall$ and $\forall \exists$.

\noindent
Each of these families is closed under
union and intersection.
But not under complementation since $\Sigma^0_{i}$ and
$\Pi^0_{i}$ are so exchanged.

\noindent A last notation:
$\Delta^0_{i}$ denotes $\Sigma^0_{i} \cap \Pi^0_{i}$.
In particular, $\Delta^0_{1}$ means c.e. and co-c.e. hence
computable.

\noindent
As for inclusion, $\Delta^0_{2}$ strictly
contains the boolean closure of $\Sigma^0_{1}$,
in particular it contains
$\Sigma^0_{1} \cup \Pi^0_{1}$.
This is why the term hierarchy is used.

\medskip
Also, we see that $K_{A}$ is quite low as a
$\Delta^0_{2}$ relation since
$\Sigma^0_{1}\wedge \Pi^0_{1}$ is the very first
level of the boolean closure of $\Sigma^0_{1}$.
%
%
\section{Random finite strings and their applications}
%
\subsection{Random versus how much random}
\label{randomfini}
\noindent {\bf Q}: \
Let's go back to the question: ``what is a random string?"

\medskip  \noindent {\bf A}: \
This is the interesting question,
but this will not be the one we shall answer.
We shall modestly consider the question:
``To what extent is $x$ random?"

\noindent
We know that $K(x) \leq |x|+O(1)$.
It is tempting to declare a string $x$ random
if $K(x) \geq |x|-O(1)$.
But what does it really mean?
The $O(1)$ hides a constant. Let's explicit it.
\begin{definition}
A string is called $c$-incompressible
(where $c\geq 0$ is any constant)
if $K(x) \geq |x|-c$.
Other strings are called $c$-compressible.

\noindent
$0$-incompressible strings are also called
incompressible.
\end{definition}

\medskip  \noindent {\bf Q}: \
Are there many $c$-incompressible strings?

\medskip  \noindent {\bf A}: \
Kolmogorov noticed that they are quite numerous.
\begin{theorem}
For each $n$ the proportion of $c$-incompressible
among strings with length $n$ is $>1- 2^{-c}$.
\end{theorem}
\noindent
For instance, if $c=4$ then, for any length $n$,
more than $90$\% of strings are $4$-incompressible.
With $c=7$ and $c=10$ we go to more than $99$\%
and $99.9$\%.

\medskip \noindent
The proof is a simple counting argument.
There are $1+2+2^2+\dots+2^{n-c-1}=2^{n-c}-1$
programs with length $< n-c$.
Every string with length $n$ which is $c$-compressible
is necessarily the output of such a program
(but some of these programs may not halt or may output
a string with length $\neq n$).
Thus, there are at most $2^{n-c}-1$ $c$-compressible
strings with length $n$,
hence at least $2^n -(2^{n-c}-1)=2^n -2^{n-c}+1$
$c$-incompressible strings with length $n$.
Whence the proportion stated in the theorem.

\medskip  \noindent {\bf Q}: \
Are $c$-incompressible strings really random?

\medskip  \noindent {\bf A}: \
Yes. Martin-L\"of, 1965 \cite{martinlof66},
formalized the notion of statistical test
and proved that incompressible strings pass all these
tests (cf. \S\ref{sub:randomseq}).
%
\subsection{Applications of random finite strings
in computer science}
\noindent {\bf Q}: \
And what is the use of incompressible strings in
computer science?

\medskip  \noindent {\bf A}: \
Roughly speaking, incompressible strings are
strings without any form of local or global
regularity.
Consideration of such objects may help
almost anytime one has to show something is
complex, for instance a lower bound for worst
or average case time/space complexity.
The accompanying key tool is Proposition
\ref{prop:K-f(x)}.

\noindent
And, indeed, incompressible strings have been
successfully used in such contexts.
An impressive compilation of such applications
can be found in Ming Li and Paul Vitanyi's book
(\cite{li-vitanyi}, chapter 6), running
through nearly 100 pages!

\medskip  \noindent {\bf Q}: \
Could you give an example?

\medskip  \noindent {\bf A}: \
Sure. The very first such application is quite
representative.
It is due to Wolfgang Paul, 1979 \cite{paul79}
and gives a quadratic lower bound
on the computation time of any one-tape
Turing machine ${\cal M}$ which recognizes palindromes.

\medskip  \noindent
Up to a linear waste of time, one can suppose that
${\cal M}$ always halts on its first cell.

\noindent
Let $n$ be even and
$x x^R = x_{1}x_{2}\ldots x_{n-1}x_{n}
x_{n}x_{n-1}\ldots x_{2}x_{1}$ be a palindrome
written on the input tape of the Turing machine
${\cal M}$.

\noindent
For each $i < n$ let $CS_{i}$ be the
crossing sequence associated to cell $i$, i.e.
the list of successive states of ${\cal M}$
when its head visits cell $i$.

\medskip  \noindent
Key fact: \
{\em string $x_{1}x_{2}\ldots x_{i}$ is
uniquely determined by $CS_{i}$.}

\noindent
I.e. $x_{1}x_{2}\ldots x_{i}$ is the sole string $y$
such that
--- relative to an ${\cal M}$-computation
on some palindrome with prefix $y$ ---,
the crossing sequence on cell $|y|$ is $CS_{i}$.

\medskip  \noindent
This can be seen as follows.
Suppose $y\neq x_{1}x_{2}\ldots x_{i}$ leads to the
same crossing sequence $CS_{i}$ on cell $|y|$
for an ${\cal M}$-computation on some palindrome
$yzz^Ry^R$.
Run ${\cal M}$ on input $yx_{i+1}\ldots x_{n}x^R$.
Consider the behaviour of ${\cal M}$
while the head is on the left part $y$.
This behaviour is exactly the same as that for
the run on input $yzz^Ry^R$ because the sole useful
information for ${\cal M}$ while scanning $y$ comes
from the crossing sequence at cell $|y|$.
In particular, ${\cal M}$ -- which halts on cell 1 ---
accepts this input $yx_{i+1}\ldots x_{n}x^R$.
But this is not a palindrome! Contradiction.

\medskip  \noindent
Observe that the way $x_{1}x_{2}\ldots x_{i}$ is
uniquely determined by $CS_{i}$ is quite complex.
But we don't care about that. It will just charge
the $O(1)$ constant in (\ref{eq:10}).

\noindent
Using Proposition \ref{prop:K-f(x)} with the
binary string associated to $CS_{i}$ which is
$c$ times longer (where $c= \lceil |Q| \rceil$,
$|Q|$ being the number of states), we see that
\begin{equation}\label{eq:10}
K(x_{1}x_{2}\ldots x_{i})
\leq c |T_{i}| + O(1)
\end{equation}
If $i \geq \frac{n}{2}$ then
$x_{1}x_{2}\ldots x_{\frac{n}{2}}$ is uniquely
determined by the pair
$(x_{1}x_{2}\ldots x_{i}, \frac{n}{2}$).
Hence also by the pair $(CS_{i}, \frac{n}{2}$).
Since the binary representation of $\frac{n}{2}$
uses $\leq \log(n)$ bits,
this pair can be encoded with
$2 c |CS_{i}| + \log(n) + 1$ bits. Thus,
\begin{equation}\label{eq:11}
K(x_{1}x_{2}\ldots x_{\frac{n}{2}})
\leq 2 c |CS_{i}| + \log(n) + O(1)
\end{equation}
Now, let's sum equations (\ref{eq:11}) for
$i=\frac{n}{2},\ldots,n$.
Observe that the sum of the lengthes of the
crossing sequences
$CS_{\frac{n}{2}},\ldots, CS_{n}$ is at most
the number $T$ of computation steps.
Therefore, this summation leads to
\begin{equation}\label{eq:12}
\frac{n}{2} K(x_{1}x_{2}\ldots x_{\frac{n}{2}})
\leq 2 c T + \frac{n}{2} \log(n) + O(\frac{n}{2})
\end{equation}
Now, consider as $x$ a string such that
$x_{1}x_{2}\ldots x_{\frac{n}{2}}$ is incompressible,
i.e.
$K(x_{1}x_{2}\ldots x_{\frac{n}{2}})
\geq \frac{n}{2}$.
Equation (\ref{eq:12}) leads to
\begin{equation}\label{eq:13}
(\frac{n}{2})^2
\leq 2 c T + \frac{n}{2} \log(n) + O(\frac{n}{2})
\end{equation}
whence $T\geq O(n^2)$. Since the input $xx^R$
has length $2n$, this proves the quadratic lower bound.
QED
%
%
\section{Prefix complexity}
%
\subsection{Self delimiting programs}
\noindent {\bf Q}: \
I heard about prefix complexity. What is it?

\medskip  \noindent {\bf A}: \
Prefix complexity is a very interesting variant
of Kolmogorov complexity which was introduced
around 1973 by Levin \cite{levin73} and,
independently, by Chaitin \cite{chaitin75}.

\noindent
The basic idea is taken from some programming
languageswhich have an explicit delimiter
to mark the end of a program.
For instance, $PASCAL$ uses ``{\em end.}".
Thus, no program can be a proper prefix of another
program.

\medskip  \noindent {\bf Q}: \
This is not true with $PROLOG$ programs:
you can always add a new clause.

\medskip  \noindent {\bf A}: \
To execute a $PROLOG$ program, you have to
write down a query.
And the end of a query is marked by a full stop.
So, it's also true for $PROLOG$.

\medskip  \noindent {\bf Q}: \
OK. However, it's not true for $C$ programs nor
$LISP$ programs.

\medskip  \noindent {\bf A}: \
Hum\ldots You are right.
%
\subsection{Chaitin-Levin prefix complexity}
\noindent {\bf A}: \
Let's say that a set $X$ of strings is prefix-free
if no string in $X$ is a proper prefix of another
string in $X$.
A programming language
$A:\{0,1\}^* \rightarrow \{0,1\}^*$
is prefix if its domain is a prefix-free set.

\medskip  \noindent {\bf Q}: \
So the programming language $PASCAL$ that you seem to be fond of
is prefix.

\medskip  \noindent {\bf A}: \
Sure, $PASCAL$ is prefix.

\medskip  \noindent {\bf A}: \
But, what's new with this special condition?

\medskip  \noindent {\bf A}: \
Kolmogorov Invariance Theorem from \S\ref{sub:invariance}
goes through with prefix programming languages, leading to
the prefix variant $H$ of $K$.
\begin{theorem}[Invariance theorem]
\label{thm:invariance-prefix}
There exists a prefix partial computable
function
$U^{prefix}:\{0,1\}^* \rightarrow\{0,1\}^*$.
such that
$K_{U^{prefix}} \leq K_{A} +O(1)$
for any prefix partial computable function
$A:\{0,1\}^* \rightarrow\{0,1\}^*$.
In other words, up to an additive constant,
$K_{U^{prefix}}$ {\em is the smallest one among the}
$K_{A}$'s.
\end{theorem}
\begin{definition}[Prefix Kolmogorov complexity]
Prefix Kolmogorov complexity $H:\{0,1\}^* \rightarrow\bbbn$
is any fixed such function $K_{U^{prefix}}$.
\end{definition}
%
\subsection{Comparing $K$ and $H$}
\noindent {\bf Q}: \
How does $H$ compare to $K$ ?

\medskip  \noindent {\bf A}: \
A simple relation is as follows:
\begin{proposition}
$K(x) - O(1) \leq H(x)
\leq K(x) + 2 \log(K(x)) + O(1)$\\
Idem with $K(\ \mid\ )$ and $H(\ \mid\ )$.
\end{proposition}
The first inequality is a mere application of the
Invariance Theorem for $K$
(since $U^{prefix}$ is a programming language).
To get the second one, we consider a
programming language $U$ such that $K=K_{U}$
and construct a prefix programming language
$U'$ as follows: the domain of $U'$ is the set of
strings of the form $pad(|p|)1p$ and
$U'(pad(|p|)1p) = U(p)$.
By very construction, the domain of $U'$ is
prefix-free.
Also, $K_{U'}(x) = K_{U} (x) + 2 \log(K_{U}(x))+1$.
An application of the Invariance Theorem for $H$
gives the second inequality of the Proposition.

\medskip  \noindent
This inequality can be improved.
A better encoding leads to\\
{\centerline{$H(x)
\leq K(x) + \log(K(x)) +2 \log\log(K(x)) + O(1)$}

\medskip  \noindent
Sharper relations have been proved by Solovay,
1975 (unpublished \cite{solovay74},
cf. also \cite{li-vitanyi} p. 211):
\begin{proposition}
$H(x) = K(x) + K(K(x)) + O(K(K(K(x))))$

\centerline{
$K(x) = H(x) - H(H(x)) - O(H(H(H(x))))$}
\end{proposition}
%
\subsection{How big is $H$ ?}
\noindent {\bf Q}: \
How big is $H$ ?

\medskip  \noindent {\bf A}: \
$K$ and $H$ behave in similar ways.
Nevertheless, there are some differences.
Essentially a logarithmic term.
\begin{proposition}
$H(x) \leq |x| + 2\log(|x|)+O(1)$
\end{proposition}
\noindent
To prove it, apply the $H$ Invariance Theorem to
the prefix function

\centerline{
$A(pad(|x|)1x)=x$.}
\noindent
Of course, it can be improved to \\
\centerline{
$H(x) \leq |x| + \log(|x|) + 2\log\log(|x|)+O(1)$}

\medskip  \noindent {\bf Q}: \
How big can be $H(x)-|x|$ ?

\medskip  \noindent {\bf A}: \
Well, to get a non trivial question,
we have to fix the length of the $x$'s.
The answer is not a simple function of $|x|$
as expected, it does use $H$ itself:\\
\centerline{
$\max_{|x|=n} (H(x)- |x|) = H(|x|)+O(1)$}

\medskip  \noindent {\bf Q}: \
How big can be $H(x)-K(x)$ ?

\medskip  \noindent {\bf A}: \
It can be quite large:\\
\centerline{
$K(x) \leq |x| - \log(|x|) \leq |x| \leq H(x)$}
happens for arbitrarily large $x$'s
(\cite{li-vitanyi} Lemma 3.5.1 p. 208).
%
\subsection{Convergence of series and the Coding Theorem}
\noindent {\bf Q}: \
What's so really special with this prefix condition?

\medskip  \noindent {\bf A}: \
The possibility to use Kraft's inequality.
This inequality tells you that if $Z$ is a
prefix-free set of strings then
$\Sigma_{p\in Z} 2^{-|p|} \leq 1$.

\noindent
Kraft's inequality is not hard to prove.
Denote $I_u$ the set of infinite strings
which admit $u$ as prefix.
Observe that \\
1) $2^{-|p|}$ is the probability of $I_u$.\\
2) If $u,v$ are prefix incomparable
then $I_u$ and $I_v$ are disjoint.\\
3) Since $Z$ is prefix, the $I_u$'s, $u\in Z$ are
pairwise disjoint and their union has probability
$\Sigma_{p\in Z} 2^{-|p|} < 1$

\medskip  \noindent
The $K_A(x)$'s are lengthes of distinct programs in
a prefix set (namely, the domain of $A$).
So, Kraft's inequality implies\\
\centerline{$\Sigma_{x\in \{0,1\}^*} 2^{-K_{A}(x)} < 1$.}

\medskip  \noindent
In fact, $H$ satisfies the following very important property,
proved by Levin \cite{levin84} (which can be seen as
another version of the Invariance Theorem for $H$):
\begin{theorem}[Coding Theorem]
Up to a multiplicative factor,
$2^{-H}$ is maximum among functions
$F:\{0,1\}^* \rightarrow {\bf R}$
such that
$\Sigma_{x\in \{0,1\}^*} F(x) < +\infty$
and which are approximable from below
(in a sense dual to that in \S \ref{sub:approx},
i.e. the set of pairs $(x,q)$ such that $q$ is
rational and $q<F(x)$ is c.e.).
\end{theorem}
%
%
\section{Random infinite sequences}
%
\subsection{Top-down approach to randomness of
             infinite sequences}
\noindent {\bf Q}: \
So, we now come to random infinite sequences.

\medskip  \noindent {\bf A}: \
It happens that there are two equivalent ways
to get a mathematical notion of random sequences.
We shall first consider the most natural one,
which is a sort of ``top-down approach".

Probability laws tell you that with probability one
such and such things happen,
i.e. that some particular set of sequences has
probability one.
A natural approach leads to consider as random
those sequences which satisfy all such laws,
i.e. belong to the associated sets
(which have probability one).

An easy way to realize this would be to declare
a sequence to be random just in case it belongs
to all sets (of sequences) having probability one
or, equivalently, to no set having probability zero.
Said otherwise, the family of random sequences
would be the intersection of all sets
having probability one, i.e. the complement
of the union of all sets having probability zero.

\noindent
{\em Unfortunately, this family is empty!}
In fact, let $r$ be any sequence:
the singleton set $\{r\}$ has probability zero
and contains $r$.

In order to maintain the idea, we have to consider
a not too big family of sets with probability one.

\medskip  \noindent {\bf Q}: \
A countable family.

\medskip  \noindent {\bf A}: \
Right. The intersection of a countable family
of set with probability one will have
probability one.
So that the set of random sequences will have
probability one, which is a much expected property.
%
\subsection{Frequency tests and
             von Mises random sequences}\label{vonmises}
\noindent {\bf A}: \
This top-down approach was pionneered by Richard von Mises
in 1919 (\cite{mises19}, \cite{mises39})
who insisted on frequency statistical tests.
He declared an infinite binary sequence
$a_{0} a_{1} a_{2}\ldots$ to be random
(he used the term {\em Kollektiv}) if the frequence
of $1$'s is ``everywhere" fairly distributed
in the following sense:

\medskip \noindent
{\em i)} Let $S_{n}$ be the number of $1$'s among
the first $n$ terms of the sequence. Then
$\lim_{n\rightarrow \infty} \frac{S_{n}}{n}=\frac{1}{2}$.

\noindent
{\em ii)} The same is true for every subsequence
$a_{n_{0}+1} a_{n_{1}+1} a_{n_{2}+1}\ldots$
where $n_{0},n_{1},n_{2}\ldots$ are the successive
integers $n$ such that
$\phi(a_{0} a_{1}\ldots a_{n})=1$ where $\phi$ is
an ``admissible" place-selection rule.

\medskip \noindent
What is an ``admissible" place-selection rule was not
definitely settled by von Mises.
Alonzo Church, 1940, proposed that admissibility
be exactly computability.

It is not difficult to prove that the family of
infinite binary sequence satisfying the above
condition has probability one for any
place-selection rule.
Taking the intersection over all computable
place-selection rules, we see that the family of
von Mises-Church random sequences has probability one.

{\em However, von Mises-Church notion of random sequence
is too large.}
There are probability laws which do not reduce to
tests with place-selection rules and are not satisfied by
all von Mises-Church random sequences.
As shown by Jean Ville, \cite{ville} 1939,
this is the case for the law of iterated logarithm.
This very important law (due to A. I. Khintchin, 1924)
expresses that with probabililty one
\begin{equation} \nonumber
\lim \sup_{n\rightarrow +\infty}
\frac{S_{n}^*}{\sqrt{2\log\log (n)}} = 1
\ \mbox{ and } \
\lim \inf_{n\rightarrow +\infty}
\frac{S_{n}^*}{\sqrt{2\log\log (n)}} = -1
\end{equation}
where
$S_{n}^*=\frac{S_{n} -\frac{n}{2}}
               {\sqrt{\frac{n}{4}}}$
(cf. William Feller's book \cite{feller},
p. 186, 204--205).

\medskip  \noindent {\bf Q}: \
Wow! What do these equations mean?

\medskip  \noindent {\bf A}: \
They are quite meaningful.
The quantities
$\frac{n}{2}$ and $\sqrt{\frac{n}{4}}$
are the expectation and standard deviation
of $S_{n}$.
So that, $S_{n}^*$ is obtained from $S_{n}$
by normalization:
$S_{n}$ and $S_{n}^*$ are linearly related as
random variables, and $S_{n}^*$'s expectation
and standard deviation are $0$ and $1$.

\noindent
Let's interpret the $\lim \sup$ equation,
the other one being similar (in fact, it can be
obtained from the first one by symmetry).

\noindent
Remember that
$\lim \sup_{n\rightarrow +\infty} f_{n}$
is obtained as follows.
Consider the sequence
$v_{n}= \sup_{m\geq n} f_{m}$.
The bigger is $n$ the smaller is the set
$\{m \ : m\geq n \ \}$.
So that the sequence $v_{n}$ decreases,
and $\lim \sup_{n\rightarrow +\infty} f_{n}$
is its limit.

\noindent
The law of iterated logarithm tells you that
with probabililty one the set
$\{n \ : \ S_{n}^* > \lambda\sqrt{2\log\log (n)} \}$
is finite in case $\lambda >1$
and infinite in case $\lambda <1$.

\medskip  \noindent {\bf Q}: \
OK.

\medskip  \noindent {\bf A}: \
More precisely, there are von Mises-Church random sequences
which satisfy
$\frac{S_{n}}{n} \geq \frac{1}{2}$ for all $n$,
a property which is easily seen to contradict the law of
iterated logarithm.

\medskip  \noindent {\bf Q}: \
So, von Mises' approach is definitely over.

\medskip  \noindent {\bf A}: \
No. Kolmogorov, 1963 \cite{kolmo63},
and Loveland, 1966 \cite{loveland66},
independently considered an extension of the notion of
place-selection rule.

\medskip  \noindent {\bf Q}: \
Kolmogorov once more\ldots

\medskip  \noindent {\bf A}: \
Indeed. Kolmogorov allows place-selection rules
giving subsequences proceeding in some new order,
i.e. mixed subsequences.
The associated notion of randomness is called
Kolmogorov stochastic randomness
(cf. \cite{kolmo-uspensky} 1987).
Since there are more conditions to satisfy,
stochastic random sequences form a subclass of
von Mises-Church random sequences.
They constitute, in fact, a proper subclass
(\cite{loveland66}).

\noindent
However, it is not known whether they satisfy all
classical probability laws.
%
\subsection{Martin-L\"of random sequences}
\label{sub:randomseq}
\noindent {\bf Q}: \
So, how to come to a successful theory of random sequences?

\medskip  \noindent {\bf A}: \
Martin-L\"of found such a theory.

\medskip  \noindent {\bf Q}: \
That was not Kolmogorov?
The same Martin-L\"of you mentioned
concerning random finite strings?

\medskip  \noindent {\bf A}: \
Yes, the same Martin-L\"of, in the very same paper
\cite{martinlof66} in 1965.
Kolmogorov looked for such a notion, but
it was Martin-L\"of, a Swedish mathematician,
who came to the pertinent idea.
At that time, he was a pupil of Kolmogorov and
studied in Moscow.
Martin-L\"of made no use of Kolmogorov
random finite string
to get the right notion of infinite random sequence.
What he did is to forget about the frequency character
of computable statistical tests (in von Mises-Church
notion of randomness) and look for what could be
the essence of general statistical tests and
probability laws.
Which he did both for finite strings
and for infinite sequences.

\medskip  \noindent {\bf Q}: \
Though intuitive, this concept is rather vague!

\medskip  \noindent {\bf A}: \
Indeed. And Martin-L\"of's analysis of what can be a
probability law is quite interesting.

\medskip
To prove a probability law amounts to prove that
a certain set $X$ of sequences has probability one.
To do this, one has to prove that the exception set
--- which is the complement
$Y=\{0,1\}^{\bbbn} \setminus X$---
has probability zero.
Now, in order to prove that
$Y \subseteq \{0,1\}^{\bbbn}$ has probability zero,
basic measure theory tells us that
one has to include $Y$ in open sets with
arbitrarily small probability.
I.e. for each $n\in \bbbn$ one must find an open set
$U_{n}\supseteq Y$
which has probability $\leq \frac{1}{2^n}$.

If things were on the real line ${\bf R}$
we would say that $U_{n}$ is a countable union of
intervals with rational endpoints.

Here, in $\{0,1\}^{\bbbn}$, $U_{n}$ is a
countable union of sets of the form
$I_{u}= u\{0,1\}^{\bbbn}$ where $u$ is a finite binary
string and $I_{u}$ is the set of infinite sequences
which extend $u$.
Well, in order to prove that $Y$ has probability zero,
for each $n\in \bbbn$ one must find a family
$(u_{n,m})_{m\in \bbbn}$ such that
$Y\subseteq \bigcup_{m} I_{u_{n,m}}$
and $Proba(\bigcup_{m} I_{u_{n,m}})\leq \frac{1}{2^n}$
for each $n\in \bbbn$.

\medskip
And now Martin-L\"of makes a crucial observation:
mathematical probability laws which we can consider
necessarily have some effective character.
And this effectiveness should reflect in the proof
as follows:\\
\centerline{\em the doubly indexed sequence
$(u_{n,m})_{{n,m\in\bbbn}}$ is computable.}

Thus, the set $\bigcup_{m} I_{u_{n,m}}$ is a
{\em computably enumerable open set} and
$\bigcap_{n} \bigcup_{m} I_{u_{n,m}}$
is a countable intersection of a
{\em computably enumerable family of open sets}.

\medskip  \noindent {\bf Q}: \
This observation has been checked for proofs of usual
probability laws?

\medskip  \noindent {\bf A}: \
Sure. Let it be the law of large numbers,
that of iterated logarithm\ldots
In fact, it's quite convincing.

\medskip  \noindent {\bf Q}: \
This open set $\bigcup_{m} I_{u_{n,m}}$ could not be
computable?

\medskip  \noindent {\bf A}: \
No. A computable set in $\{0,1\}^{\bbbn}$ is always a
{\em finite} union of $I_{u}$'s.

\medskip  \noindent {\bf Q}: \
Why?

\medskip  \noindent {\bf A}: \
What does it mean that $Z\subseteq \{0,1\}^{\bbbn}$
is computable?
That there is some Turing machine such that,
if you write an infinite sequence $\alpha$
on the input tape then after finitely many steps,
the machine tells you if $\alpha$ is in $Z$ or not.
When it does answer, the machine has read but
a finite prefix $u$ of $\alpha$, so that
it gives the same answer if $\alpha$ is replaced by
any $\beta\in I_{u}$.
In fact, an application of K\"onig's lemma
(which we shall not detail) shows that
we can bound the length of such a prefix $u$.
Whence the fact that $Z$ is a finite union of $I_{u}$'s.

\medskip  \noindent {\bf Q}: \
OK. So, we shall take as random sequences those
sequences which are outside any set which is
a countable intersection of a
computably enumerable family of open sets
and has probability zero.

\medskip  \noindent {\bf A}: \
This would be too much. Remember,
$Proba(\bigcup_{m} I_{u_{n,m}})\leq \frac{1}{2^n}$.
Thus, the way the probability of
$\bigcup_{m} I_{u_{n,m}}$
tends to $0$ is computably controlled.

\medskip
So, here is Martin-L\"of's definition:
\begin{definition}
A set of infinite binary sequences is constructively of probability
zero if it is included in $\bigcap_{n} \bigcup_{m} I_{u_{n,m}}$
where $(m,n)\mapsto u_{n,m}$ is a partial computable
function $\bbbn^2 \rightarrow \{0,1\}^*$ such that
$Proba(\bigcup_{m} I_{u_{n,m}})\leq \frac{1}{2^n}$
for all $n$.
\end{definition}
And now comes a very surprising theorem
(Martin-L\"of, \cite{martinlof66}, 1966):
\begin{theorem}
There is a largest set of sequences (for the inclusion ordering)
which is constructively of probability zero.
\end{theorem}
\noindent {\bf Q}: \
Largest? up to what?

\medskip  \noindent {\bf A}: \
Up to nothing.
Really largest set:
it is constructively of probability zero
and contains any other set constructively of
probability zero.

\medskip  \noindent {\bf Q}: \
How is it possible?

\medskip  \noindent {\bf A}: \
Via a diagonalization argument.
The construction has some technicalities
but we can sketch the ideas.
 From the well-known existence of universal c.e. sets,
we get a computable enumeration
$((O_{i,j})_{i})_{j}$ of  sequences of c.e. open sets.
A slight transformation allows to satisfy the
inequality $Proba(O_{i})\leq \frac{1}{2^i}$.
Now, set $U_{j}=\bigcup_{e} O_{e,e+j+1}$
(here lies the diagonalization!)
Clearly,
$Proba(O_{j})\leq \sum_{e} \frac{1}{2^{e+j+1}}
=\frac{1}{2^j}$,
so that $\bigcap_{j} U_{j}$ is
constructively of probability zero.
Also, $U_{j}\supseteq O_{i,j}$ for all $j\geq i$
whence
$(\bigcap_{j} U_{j})\supseteq
(\bigcap_{j} O_{i,j})$.

\medskip \noindent {\bf Q}: \
So, Martin-L\"of random sequences are exactly those lying
in this largest set.

\noindent {\bf A}: \
Yes. And all theorems in probability theory can be
strengthened by replacing ``with probability one" by
``for all Martin-L\"of random sequences"
%
\subsection{Bottom-up approach to randomness of infinite sequences:
            Martin-L\"of's Large oscillations theorem}
\noindent {\bf Q}: \
So, now, what is the bottom-up approach?

\medskip \noindent {\bf A}: \
This approach looks at the asymptotic
algorithmic complexity of the prefixes
of the infinite binary sequence
$a_{0} a_{1} a_{2}\ldots$, namely
the $K(a_{0}\ldots a_{n})$'s.

The next theorem is the first significative result
relevant to this approach.
Point 2 is due to Albert Meyer and
Donald Loveland, 1969 \cite{loveland69} p. 525.
Points 3,4 are due to Gregory Chaitin,
1976 \cite{chaitin76}.
(Cf. also \cite{li-vitanyi} 2.3.4 p.124).
\begin{theorem}\label{thm:meyer}
The following conditions are equivalent:

1) \ $a_{0} a_{1} a_{2}\ldots$ is computable

2) \ $K(a_{0}\ldots a_{n}\mid n)=O(1)$.

3) \ $|K(a_{0}\ldots a_{n})-K(n)|\leq O(1)$.

4) \ $|K(a_{0}\ldots a_{n})-\log(n)|\leq O(1)$.
\end{theorem}

\noindent {\bf Q}: \
Nice results.
Let me tell what I see.
We know that $K(x)\leq |x|+O(1)$.
Well, if we have the equality,
$K(a_{0}\ldots a_{n}) = n - O(1)$,
i.e. if maximum complexity occurs for all prefixes,
then the sequence $a_{0} a_{1} a_{2}\ldots$ should be random!
Is it indeed the case?

\medskip \noindent {\bf A}: \
That's a very tempting idea. And Kolmogorov
had also looked for such a characterization.
Unfortunately, as Martin-L\"of proved around 1965
(1966, \cite{martinlof71}),
{\em there is no such sequence}!
It is a particular case of a more general result
(just set f(n)=constant).
\begin{theorem}[Large oscillations, \cite{martinlof71}]
      \label{thm:osc}
Let $f:\bbbn \rightarrow\bbbn$ be a computable function
such that $\Sigma_{n\in\bbbn}  2^{-f(n)}=+\infty$.
Then, for {\em every} binary sequence
$a_{0} a_{1} a_{2}\ldots$
there are infinitely many $n$'s such that
$K(a_{0}\ldots a_{n}\mid n) < n - f(n)$.
\end{theorem}

\noindent {\bf Q}: \
So, the bottom-up approach completely fails as concerns
a characterization of random sequences.
Hum\ldots
But it does succeed as concerns computable sequences,
which were already fairly well characterized.
Funny!

\medskip \noindent {\bf A}: \
It's however possible to sandwich the set of
Martin-L\"of random sequences between two sets of
probability one defined in terms of the $K$
complexity of prefixes.
\begin{theorem}[\cite{martinlof71}]
Let $f:\bbbn \rightarrow\bbbn$ be computable
such that the series $\Sigma 2^{-f(n)}$ is
computably convergent.
Set
\begin{eqnarray*}
X & = & \{ a_{0} a_{1}\ldots  : \
K(a_{0}\ldots a_{n}\mid n) \geq n - O(1)
\mbox{ for infinitely many } n \mbox{'s.} \}\\
Y_{f} & = & \{ a_{0} a_{1}\ldots  : \
K(a_{0}\ldots a_{n}\mid n) \geq n - f(n)
\mbox{ for all but finitely many } n \mbox{'s.} \}
\end{eqnarray*}
Denote $ML$ the set of Martin-L\"of random sequences.
Then $X$ and $Y_{f}$ have probability one and
$X \subset ML \subset Y_{f}$.
\end{theorem}
\noindent
NB: Proper inclusions have been proved by Peter Schnorr,
1971 \cite{schnorr71}
(see also \cite{li-vitanyi} 2.5.15 p.154).

\medskip
Let's illustrate this theorem on an easy and spectacular
corollary which uses the fact that
$2^{-2\log(n)}=\frac{1}{n^2}$ and that
the series $\Sigma\frac{1}{n^2}$ is
computably convergent:
{\em if $K(a_{0}\ldots a_{n}\mid n) \geq n -c$
for infinitely many $n$'s then
$K(a_{0}\ldots a_{n}\mid n) \geq n - 2\log(n)$ for all but
finitely many $n$'s.}
%
%

%
%

%
\subsection{Bottom-up approach with prefix complexity}
  \noindent {\bf Q}: \
What about considering prefix Kolmogorov
complexity?

\medskip  \noindent {\bf A}: \
Kolmogorov's original idea does work with prefix Kolmogorov
complexity.
This has been proved by Claus Peter Schnorr
(1974, unpublished, cf. \cite{chaitin75} Remark p. 106,
and \cite{chaitin87} p. 135-137 for a proof).

\noindent
Robert M. Solovay, 1974 (unpublished \cite{solovay74})
strengthened Schnorr's result
(cf. \cite{chaitin87} p. 137-139).
\begin{theorem}\label{martinlof-H}
The following conditions are equivalent:

\noindent
1) \ $a_{0} a_{1} a_{2}\ldots$ is Martin-L\"{o}f ramdom
random

\noindent
2) \ $H(a_{0} \ldots a_{n}) \geq n-O(1)$ for all $n$.

\noindent
3) \ $\lim_{n \rightarrow +\infty}
(H(a_{0} \ldots a_{n}) - n) = +\infty$.

\noindent
4) \ For any c.e. sequence $(A_{i})_{i}$ of open subsets
of $\{0,1\}^{\bbbn}$ if
$\Sigma_{i}Proba(A_{i}) < +\infty$ then $a_{0} a_{1} a_{2}\ldots$
belongs to finitely many $A_{i}$'s.
\end{theorem}
\noindent
These equivalences stress the robustness of the notion
of Martin-L\"{o}f ramdom sequence.
%
\subsection{Top-down/Bottom-up approaches: a sum up}
\noindent {\bf Q}: \
I get somewhat confused with these two approaches.
Could you sum up.

\medskip  \noindent {\bf A}: \
The top-down and bottom-up approaches
both work and lead to the very same class of
random sequences.

\noindent
Kolmogorov looked at the bottom-up approach
from the very beginning in 1964.
But nothing was possible with
the original Kolmogorov complexity,
Levin-Chaitin's variant $H$ was needed.

\medskip  \noindent {\bf Q}: \
Ten years later\ldots

\medskip  \noindent {\bf A}: \
As for the top-down approach, it was pionneered
by von Mises since 1919 and made successful
by Martin-L\"of in 1965.
Martin-L\"of had to give up von Mises frequency tests.
However, Kolmogorov was much interested by these
frequency tests (\cite{kolmo63}),
and he refined them in a very clever way
with the purpose to recover
Martin-L\"of randomness, which
lead him to the notion of Kolmogorov
stochastic randomness.
Unfortunately, up to now, we only know that

\centerline{Martin-L\"of random
$\Rightarrow$
stochastic random
$\Rightarrow$
von Mises-Church random.}
\noindent
The second implication is known to be strict
but not the first one.
Would it be an equivalence, this would give a
quite vivid characterization of random
sequences via much concrete tests.
%
\subsection{Randomness with other probability
             distributions}
\noindent {\bf Q}: \
All this is relative to the uniform probability
distribution. Can it be extended to arbitrary
probability distributions?

\medskip  \noindent {\bf A}: \
Not arbitrary probability distributions, but
computable Borel ones: those distributions $P$
such that the sequence of reals
$(P(I_{u}))_{u\in\{0,1\}^*}$
(where $I_{u}$ is the set of infinite sequences
which extend $u$)
is computable, i.e. there is a computable function
$f:\{0,1\}^*\times \bbbn \rightarrow {\bf Q}$ such that

\centerline{$|P(I_{u}) - f(u,n)| \leq \frac{1}{2^n}$.}
\noindent
Martin-L\"of's definition of random sequences extends
trivially.
As for characterizations with variants of Kolmogorov
complexity, one has to replace the length of a finite
string $u$ by the quantity $-log(P(I_{u}))$.
%
%
\subsection{Chaitin's real $\Omega$}
\noindent {\bf Q}: \
I read a lot of things about Chaitin's real $\Omega$.

\medskip  \noindent {\bf A}: \
Gregory Chaitin, 1987 \cite{chaitin87a},
explicited a spectacular random real and
made it very popular.
\\
Consider a universal prefix partial recursive function $U$
and let $\Omega$ be the Lebesgue measure of the set
$$
\{\alpha\in\cantor \mid
\exists n\ U(\alpha\segment n)\mbox{ is defined}\}
$$

\noindent {\bf Q}: \
Seems to be an avatar of the halting problem.

\medskip  \noindent {\bf A}: \
Indeed. It is the probability that, on an infinite input,
the machine which computes $U$ halts in finite time
(hence after reading a finite prefix of its input).
\begin{equation}
\Omega = \Sigma \{ 2^{-|p|} \mid \mbox{$U$ halts on input $p$}\}
\end{equation}
\begin{theorem}
The binary expansion of $\Omega$ is Martin-L\"{o}f ramdom.
\end{theorem}

\noindent {\bf Q}: \
How does one prove that $\Omega$ is ramdom?

\medskip  \noindent {\bf A}: \
$U$ has prefix-free domain, hence
$\Omega=\Sigma \{ 2^{-|p|} \mid p\in domain(U)\}<1$.
Any halting program with length $n$ contributes for exactly
$2^{-n}$ to $\Omega$.
Thus, if you know the first $k$ digits of $\Omega$
then you know the number of halting programs with length $\leq k$.
From this number, by dovetailing, you can get the list
of the halting programs with length $\leq k$
(cf. \S \ref{sub:approx}, \ref{sub:dovetailing}).
Having these programs, you can get the first string $u$
which is not the output of such a program.
Clearly, $H(u)>k$.
Now, $u$ is computably obtained from the
first $k$ digits of $\Omega$, so that
by Proposition \ref{prop:K-f(x)} we have
$H(u) \leq H(\omega_{0} \ldots \omega_{k}) +O(1)$.
Whence
$H(\omega_{0} \ldots \omega_{k}) \geq k +O(1)$,
which is condition 2 of Theorem \ref{martinlof-H}
(Schnorr condition).
This proves that the binary expansion of
$\Omega$ is a Martin-L\"{o}f ramdom sequence.

\medskip  \noindent {\bf Q}: \
$\Omega$ seems to depend on the universal machine.

\medskip  \noindent {\bf A}: \
Sure. We can speak of the class of Chaitin $\Omega$ reals:
those reals which express the halting probability
of some universal prefix programming language.

\noindent
Cristian Calude \& Peter Hertling \&
Bakhadyr Khoussainov \& Yongge Wang, 1998
\cite{calude98}
(cf. also Anton{\`\i}n Ku\u{c}era \& Theodore Slaman,
2000 \cite{kucera-slaman})
proved a very beautiful result: $r$ is a
Chaitin $\Omega$ real if and only if
(the binary development of) $r$ is Martin-L\"of random
and $r$ computably enumerable from below
(i.e. the set of rational numbers $< r$ is c.e.).

\medskip  \noindent {\bf Q}: \
I read that this real has incredible properties.

\medskip  \noindent {\bf A}: \
This real has a very simple and appealing definition.
Moreover, as we just noticed, there is a simple way
to get all size $n$ halting programs from its $n$ first
digits. This leads to many consequences due to the following
fact:
any $\Sigma_1^0$ statement of the form
$\exists \vec x \Phi(\vec x)$
(where $\Phi$ is a computable relation)
is equivalent to a statement insuring that
a certain program halts,
and this program is about the same size as the statement.
Now, deciding the truth of $\Sigma_1^0$ statements is the same
as deciding that of $\Pi_1^0$ statements.

\noindent
And significant $\Pi_1^0$ statements abound!
Like Fermat's last theorem (which is now Wiles' theorem)
or consistency statements.
This is why Chaitin says $\Omega$ is the ``Wisdom real".

Other properties of $\Omega$ are common to all reals
which have Martin-L\"{o}f ramdom binary expansions.
For instance, transcendance and the fact
that any theory can give us but finitely many digits.

\noindent
Hum\ldots About that last point, using
Kleene's recursion theorem, Robert Solovay, 1999
\cite{solovay99}, proved that there are particular
Chaitin $\Omega$ reals about which a given theory
can not predict {\em any} single bit!
%
\subsection{Non computable invariance}
\noindent {\bf Q}: \
In some sense, Martin-L\"{o}f ramdomness is a part of
recursion theory. Do random sequences form a
Turing degree or a family of Turing degrees?

\medskip  \noindent {\bf A}: \
Oh, no! Randomness is definitely
not computable invariant.
It's in fact a very fragile notion:
quite insignificant modifications destroy randomness.
This makes objects like $\Omega$ so special.

\noindent
Let's illustrate this point on an example.
Suppose you transform a
random sequence $a_0 a_1 a_2 a_3 \ldots$
into $a_0 0 a_1 0 a_2 0 a_3 0 \ldots$
The sequence you obtain has the same Turing
degree as the original one,
but it is no more random since its
digits with odd ranks are all $0$.
A random sequence has to be random everywhere.
Hum \ldots for Martin-L\"{o}f random reals,
I should rather say "every c.e. where".

\medskip  \noindent {\bf Q}: \
Everywhat?

\medskip  \noindent {\bf A}: \
"Every c.e. where". I mean that if $f$ is a
computable function from $\bbbn$ into $\bbbn$
(in other words, a computable enumeration of
an c.e. set)
then the sequence of digits with ranks
$f(0),f(1),f(2),\ldots$ of a Martin-L\"{o}f
random sequence has to be Martin-L\"{o}f random.
In fact, you recognize here an extraction
process \`a la von Mises for which a random sequence
should give another random sequence.

\medskip  \noindent {\bf Q}: \
OK. What about many-one degrees?

\medskip  \noindent {\bf A}: \
Same. Let's represent a binary infinite sequence
$\alpha$ by the set $X_{\alpha}$ of positions of
digits $1$ in $\alpha$.
Then,\\
\centerline{
$n\in X_{a_0 a_1 a_2\ldots} \Leftrightarrow
  2n \in X_{a_0 0 a_1 0 a_2 0 \ldots}$}
Also, let $\varphi(2n)=n$ and $\varphi(2n+1)=k$ where
$k$ is some fixed rank such that $a_k=0$, then\\
\centerline{
$n\in X_{a_0 0 a_1 0 a_2 0 ldots} \Leftrightarrow
\varphi(n) \in X_{a_0 a_1 a_2\ldots}$}
These two equivalences prove that $X_{a_0 a_1 a_2\ldots}$
and $X_{a_0 0 a_1 0 a_2 0\ldots}$ are many-one equivalent.
%
%
\section{More randomness}
\begin{quote}
{\em There are more things in heaven and earth, Horatio,\\
Than are dreamt of in your philosophy.}
\end{quote}
\begin{flushright}
{\em Hamlet}, William Shakespeare
\end{flushright}
%
\subsection{Beyond c.e.: oracles and infinite computations}
\noindent {\bf Q}: \
Are there other random reals than
Chaitin $\Omega$ reals ?

\medskip  \noindent {\bf A}: \
Sure. Just replace in Martin-L\"{o}f's definition
the computable enumerability condition by a more complex one.
For instance, you can consider $\Sigma_2^0$ sets,
which amounts to computable enumerability with oracle $\emptyset'$
(the set which encodes the halting problem for Turing machines).

\medskip  \noindent {\bf Q}: \
Wait, wait. Just a minute ago, you said
that for all classical probability laws,
c.e. open sets, i.e. $\Sigma_1^0$ sets,
are the ones which come in
when proving that the exception set
to the law has probability zero.
So, what could be the use of such
generalizations?

\medskip  \noindent {\bf A}: \
Clearly, the more random sequences you have
which satisfy classical probability laws,
the more you strengthen these theorems as
we said earlier.
In this sense, it is better to stick to
Martin-L\"{o}f's definition.
But you can also want to consider random
sequences as worst objects to use in some context.
Depending on that context, you can be lead
to ask for much complex randomness conditions.

Also, you can have some very natural objects
much alike Chaitin real $\Omega$ which can be more complex.

\medskip  \noindent {\bf Q}: \
Be kind, give an example!

\medskip  \noindent {\bf A}: \
In a recent paper, 2001 \cite{alpha},
Ver\'onica Becher \& Chaitin \& Sergio Daicz consider
the probability that a prefix universal programming
language produces a finite output,
though possibly running indefinitely.
They prove that this probability is an $\Omega'$ Chaitin real,
i.e. its binary expansion is $\emptyset'$-random.
Becher \& Chaitin, 2002 \cite{becher-chaitin},
consider the probability for the output to represent a
cofinite set of integers, relatively to some coding of
sets of integers by sequences.
They prove it to be an $\Omega''$ Chaitin real.
\\
Such reals are as much appealing and remarkable as Chaitin's real
$\Omega$ and also they are logically more complex.
%
%
\subsection{Far beyond: Solovay random reals in set theory}
\noindent {\bf Q}: \
I heard about Solovay random reals in set theory.
Has it anything to do with Martin-L\"{o}f
random reals?

\medskip  \noindent {\bf A}: \
Hum... These notions come from very different
contexts. But well, there is a relation:
proper inclusion. Every Solovay random real
is Martin-L\"{o}f random. The converse being
far from true.
In fact, these notions of randomness are two
extreme notions of randomness.
Martin-L\"{o}f randomness is the weakest
condition whereas Solovay randomness is really
the strongest one. So big indeed that for Solovay
reals you need to work in set theory, not merely
in recursion theory, and even worse, you have to
consider two models of set theory,
say $M_{1}$ and an inner submodel $M_{2}$ with the same
ordinals\ldots

\medskip  \noindent {\bf Q}: \
You mean transfinite ordinals?

\medskip  \noindent {\bf A}: \
Yes,
$0,1,2,3,\ldots,\omega,
\omega+1,\omega+2,\ldots,\omega+\omega$
(which is $\omega.2$) an so on:
$\omega.3,\ldots, \omega.\omega$
(which is $\omega^2$)
$,\ldots, \omega^3,\ldots, \omega^\omega, \ldots$

\noindent In a model of set theory, you have reals and may
consider Borel sets, i.e. sets obtained from
rational intervals via iterated countable unions
and countable intersections.

\noindent Thus, you have reals in $M_{1}$ and reals in $M_{2}$
and every $M_{2}$ real is also in $M_{1}$.
You also have Borel sets defined in $M_{2}$. And
to each such Borel set $X_2$ corresponds a Borel
set $X_1$ in $M_{1}$ with the same definition (well,
some work is necessary to get a precise meaning,
but it's somewhat intuitive).
One can show that $X_2 \subseteq X_1$ and that
$X_1,X_2$ have the very same measure, which is
necessarily a real in $M_{2}$.
Such a Borel set $X_1$ in $M_{1}$ will be called a
$M_{2}$-coded Borel set.

\noindent Now, a real $r$ in $M_1$ is Solovay random
over $M_{2}$ if it lies in no measure zero $M_{2}$-coded
Borel set of $M_{1}$.
Such a real $r$ can not lie in the inner model
$M_{2}$ because $\{r\}$ is a measure zero Borel set
and if $r$ were in $M_{2}$ then $\{r\}$
would be $M_{2}$-coded and $r$ should be outside it,
a contradiction.

\noindent In case $M_{1}$ is big enough relative to $M_{2}$
it can contain reals which are Solovay random over $M_{2}$.
It's a rather tough subject, but you see:

\noindent --- Martin-L\"{o}f random reals are reals outside
all c.e. $G_{\delta}$ sets
(i.e. intersection of an c.e. sequence of open sets)
constructively of measure zero.
In other words, outside a very smooth countable family
of Borel sets.
Such Borel sets are, in fact, coded in any inner
submodel of set theory.

\noindent --- Solovay random reals over a submodel of
set theory are reals outside every measure zero Borel set
coded in that submodel.
Thus Solovay random reals can not be in the inner submodel.
They may or may not exist, depending on how big is $M_1$
relative to $M_2$.

\medskip  \noindent {\bf Q}: \
What a strange theory. What about the motivations?

\medskip  \noindent {\bf A}: \
Solovay introduced random reals
in set theory at the pionneering time of
independence results in set theory,
using the method of forcing invented by Paul J. Cohen.
That was in the 60's.
He used them to get a model of set theory in which
every set of reals is Lebesgue measurable \cite{solovay}.

\medskip  \noindent {\bf Q}: \
Wow! it's getting late.

\medskip  \noindent {\bf A}: \
  Hope you are not exhausted.

\medskip  \noindent {\bf Q}: \
I really enjoyed talking with you on such a topic.

\vspace{6mm}
\noindent {\em Note.}
The best reference to the subject are
\begin{itemize}
\item
Li \& Vitanyi's book \cite{li-vitanyi}
\item
Downey \& Hirschfeldt's book \cite{downeybook}
\end{itemize}
Caution: In these two books, $C,K$ denote what is here
-- and in many papers -- denoted $K,H$.
\\
Among other very useful references:
\cite{calude}, \cite{delahaye}, \cite{gacs93},
\cite{shen} and \cite{uspensky-semenov-shen}.
\medskip\\
Gregory Chaitin's papers are available on his home page.
\nocite{*}
\bibliographystyle{plain}


\end{document}